\documentclass{conm-p-l}

\newtheorem{theorem}{Theorem}[section]
\newtheorem{lemma}[theorem]{Lemma}
\newtheorem{proposition}[theorem]{Proposition}
\newtheorem{corollary}[theorem]{Corollary}
\newtheorem{conjecture}[theorem]{Conjecture}

\theoremstyle{definition}
\newtheorem{definition}[theorem]{Definition}

\theoremstyle{remark}
\newtheorem{remark}[theorem]{Remark}

\numberwithin{equation}{section}

\def \End{{\rm End}}

\def \Hom{{\rm Hom}}
\def \Map{{\rm Map}}
\def \wt {{\rm wt}}

\begin{document}

\title{On ${\mathbb Z}$-graded associative algebras and their
${\mathbb N}$-graded modules}

\author{Haisheng Li}
\address{Department of Mathematical Sciences, 
Rutgers University, Camden, NJ 08102}
\email{hli@crab.rutgers.edu}
\thanks{H. Li was supported in part by NSF Grant DMS-9616630.}

\author{Shuqin Wang}
\address{Department of Mathematical Sciences,
Harbin Institute of Technology at Weihai,
Weihai, Shandong 264209, P. R. China}



\begin{abstract}
Let $A$ be a ${\mathbb Z}$-graded associative algebra
and let $\rho$ be an irreducible  ${\mathbb N}$-graded representation of $A$ 
on $W$ with finite-dimensional homogeneous subspaces. Then it is proved 
that $\rho(\tilde{A})=gl_{J}(W)$, where $\tilde{A}$ 
is the completion of $A$ with respect to a certain topology
and $gl_{J}(W)$ is the subalgebra of $\End W$, generated by homogeneous
endomorphisms. It is also proved that an ${\mathbb N}$-graded vector space
$W$ with finite-dimensional homogeneous spaces
is the only continuous irreducible ${\mathbb N}$-graded
$gl_{J}(W)$-module up to equivalence, where $gl_{J}(W)$ 
is considered as a topological algebra in a certain natural way,
and that any continuous 
${\mathbb N}$-graded $gl_{J}(W)$-module
is a direct sum of some copies of $W$. A duality for certain 
subalgebras of $gl_{J}(W)$ is also obtained.
\end{abstract}

\maketitle

\section{Introduction}
This paper was motivated by a paper of Malikov [M] and
by the theory of vertex operator algebras.
Malikov considered the following situation:
Let $A_{1}$ and $A_{2}$ be associative algebras with unit 
over a field ${\mathbb F}$ and let
$\rho: A_{1}\otimes A_{2} \rightarrow \End W$
be an irreducible representation of $A_{1}\otimes A_{2}$ on
$W$. Naturally identify
$A_{1}$ and $A_{2}$ with subalgebras $A_{1}\otimes {\mathbb F}$ and
${\mathbb F}\otimes A_{2}$ of $A_{1}\otimes A_{2}$, respectively.
Set
$$\rho(A_{i})^{\wedge}=\Hom_{A_{i}}(W,W)\;\;\;\;(i=1,2).$$
Clearly, $\rho(A_{1})\subset \rho(A_{2})^{\wedge}, \;\;
\rho(A_{2})\subset \rho(A_{1})^{\wedge}$.
If $W$ is {\em finite-dimensional} and ${\mathbb F}$ is algebraically closed, 
the classical Brauer theorem (a special case) asserts that 
if $W$ is a completely reducible
$A_{1}$-module, then $W$ is also a completely reducible
$A_{2}$-module and the following duality holds:
$$\rho(A_{1})=\rho(A_{2})^{\wedge}, 
\;\;\rho(A_{2})=\rho(A_{1})^{\wedge}.$$
A particular case of the Brauer theorem with $A_{2}={\mathbb F}$ 
gives the Burnside theorem which asserts that $\rho(A)=\End W$. 

As pointed out in [M], if $W$ is {\em infinite-dimensional},
even the Burnside theorem is not true any more.
Nevertheless, Malikov established certain infinite-dimensional
analogues of the two classical theorems for $A$ being a 
${\mathbb Z}$-graded associative algebra of a certain type. 
As we shall explain below, on one hand, he considered
a certain (smaller) subalgebra $gl_{J}(W)$ of $\End W$
and on the other hand, he considered 
the completion $\tilde{A}$ (a ``bigger'' algebra) of a 
${\mathbb Z}$-graded algebra $A$ with respect to a certain topology.

Let $W=\coprod_{n\ge 0}W(n)$ be an ${\mathbb N}$-graded vector space  
with all $W(n)$ being finite-dimensional.
For $m\in {\mathbb Z}$, denote by
$(\End W)_{m}$ the space of all homogeneous endomorphisms of $W$
of degree $m$. Set
$$gl_{J}(W)=\coprod_{m\in {\mathbb Z}}(\End W)_{m}.$$
Then $gl_{J}(W)$ is a subalgebra of $\End W$ and it is a ${\mathbb Z}$-graded
algebra itself. 
Now let $A=\coprod_{n}A_{n}$ be any ${\mathbb Z}$-graded associative algebra.
For $n,k\in {\mathbb Z}$, set $A_{n,k}=\sum_{m>k}A_{n+m}A_{-m}$.
Endow $A_{n}$ with a topology with $a+A_{n,k}$ for $a\in A_{n},\;k\in {\mathbb Z}$
as a base of open sets (cf. [FZ], [M]).
In this way, we have a topology on $A$ and
$A$ becomes a topological algebra.
Let $\tilde{A}_{n}$ be the completion of $A_{n}$ and set
$\tilde{A}=\coprod_{n}\tilde{A}_{n}$. 
Then $\tilde{A}$ is a topological algebra.
With these notions, any ${\mathbb N}$-graded representation $\rho$ 
of $A$ on $W$ naturally extends to a representation, 
denoted by $\tilde{\rho}$,
of $\tilde{A}$ on $W$ and $\tilde{\rho}(\tilde{A})\subset gl_{J}(W)$.
Malikov proved that if $A$ is what he called a ``${\mathbb Z}$-graded algebra 
with involution,'' then $\tilde{\rho}(\tilde{A})=gl_{J}(W)$. 
Using this, Malikov easily established an analogue of the Brauer theorem
(a special case).

{}From the definition in [M], a ${\mathbb Z}$-graded algebra $A$ with involution
has the main features of the universal enveloping algebras
of certain Lie algebras such as Kac-Moody algebras, 
Heisenberg algebras and the Virasoro algebra.
Specifically, there exist a
unit, a counit, a triangular decomposition
$A=A^{+}\otimes A^{0}\otimes A^{-}$, where $A^{\pm}$ is a subalgebra of 
$A_{\pm}=\coprod_{n\ge 1}A_{\pm n}$ and $A^{0}$ is a commutative subalgebra 
of $A_{0}$, and an involution $\omega$ such that 
$\omega|_{A^{0}}=1, \;\omega (A^{\pm})=A^{\mp}$.
Lie algebras such as affine Kac-Moody algebras, 
Heisenberg algebras and the Virasoro algebra are 
known to be important sources of
vertex operator algebras (cf. [DLe], [FLM], [FZ], [Li]).
On the other hand, certain ${\mathbb Z}$-graded associative algebras
also naturally come out
in the study of vertex operator algebras (cf. [FZ], [DLin], [KL]). 
For example, to each vertex operator algebra $V$
Frenkel and Zhu [FZ] associated a ${\mathbb Z}$-graded topological associative algebra
$U(V)$
such that the category of ${\mathbb N}$-graded weak $V$-modules is equivalent to
the category of continuous ${\mathbb N}$-graded $U(V)$-modules.
For {\em certain} vertex operator algebras $V$, one may prove that
$U(V)$ is a ${\mathbb Z}$-graded algebra with involution,
so that Malikov's result can be applied.
{\em Instead}, our main goal here is to prove Malikov's analogue of
the Burnside theorem for an {\em arbitrary} ${\mathbb Z}$-graded algebra $A$
and give an analogue (more general than Malikov's) of the Brauer's theorem.
It is our belief that Frenkel and Zhu's universal enveloping algebra 
$U(V)$ is an appropriate device in the study of
some conjectured duality in vertex operator algebra theory.
This is also our main motivation of this paper.

The main results and the organization of this paper are described 
as follows: In Section 2 we review
a certain completion of ${\mathbb Z}$-graded associative algebras,
define the notion of continuous module and define the topological algebra 
$gl_{J}(W)$. 
In Section 3, motivated by the notion of rationality 
(cf. [Z], [DLM2]) in vertex operator algebra theory, we prove
that $W$ is the only continuous irreducible
${\mathbb N}$-graded module for $gl_{J}(W)$
up to equivalence and that any continuous
${\mathbb N}$-graded module is completely reducible. 
This shows that $gl_{J}(W)$ resembles 
(finite-dimensional) matrix algebras. 
In Section 4 we prove Malikov's analogue of the Burnside theorem
for an arbitrary ${\mathbb Z}$-graded associative algebra.
In Section 5 we  give a certain duality  as an analogue of
the Brauer's theorem. In Section 6 we give an application
in vertex operator algebra theory.

One of us, H. Li, would like to thank Martin Karel for many
useful discussions and collaboration in [KL].

\section{A certain completion of a ${\mathbb Z}$-graded associative algebra}
Throughout this section,
$W=\coprod_{n\in {\mathbb N}}W(n)$ is an ${\mathbb N}$-graded vector space
such that $\dim W(n)<\infty$ for $n\in {\mathbb N}$,
$A=\coprod_{n\in {\mathbb Z}}A_{n}$ represents a ${\mathbb Z}$-graded
 associative algebra and $B=\coprod_{n\in {\mathbb Z}}B_{n}$ represents
a ${\mathbb Z}$-graded topological associative algebra.

An element $f$ of $\End W$ is said to be {\em homogeneous of degree
$k$} if
\begin{eqnarray}
fW(n)\subset W(n+k)\;\;\;\mbox{ for }n\in {\mathbb Z}.
\end{eqnarray}
Let $gl_{J}(W)_{k}$ be the space of all homogeneous endomorphisms of $W$ of 
degree $k$ and set
\begin{eqnarray}\label{eglsum}
gl_{J}(W)=\coprod_{k\in {\mathbb Z}}gl_{J}(W)_{k}.
\end{eqnarray}
Clearly, $gl_{J}(W)$ is an (associative)  subalgebra of $\End W$
and it is a ${\mathbb Z}$-graded unital algebra itself 
with respect to the grading given in (\ref{eglsum}).

For $n\in {\mathbb Z}$, let $p_{n}$ be the projection map of $W$ onto $W(n)$.
Then $f\in gl_{J}(W)$ if and only if there is a nonnegative integer $r$ 
such that
\begin{eqnarray}
fW(n)\subset \oplus_{m=n-r}^{n+r}W(m)\;\;\;\mbox{ for }n\in {\mathbb Z},
\end{eqnarray}
or equivalently,
\begin{eqnarray}\label{e2.4}
p_{m}fp_{n}=0\;\;\;\mbox{ if } |m-n|>r.
\end{eqnarray}
In literatures (cf. [KP], [M]), (\ref{e2.4}) was
commonly used to define $gl_{J}(W)$.

Since $\End W=\prod_{m\in {\mathbb Z}}\Hom(W(m),W)$, we have
\begin{eqnarray}\label{e2.5}
gl_{J}(W)_{n}=\prod_{m\in {\mathbb Z}}\Hom(W(m),W(m+n)).
\end{eqnarray}
In view of this, we may (and we shall) consider $gl_{J}(W)$ as a 
completion. 
For $n\in {\mathbb Z}$, we define
\begin{eqnarray}
(\End'W)_{n}=\coprod_{m\in {\mathbb Z}}\Hom(W(m),W(m+n)).
\end{eqnarray}
Then set
\begin{eqnarray}\label{eendf1}
\End' W=\coprod_{n\in {\mathbb Z}}(\End'W)_{n}.
\end{eqnarray}
That is,
\begin{eqnarray}
\End' W=\coprod_{m,n\in {\mathbb Z}}\Hom(W(m),W(n)).
\end{eqnarray}
In terms of the projection maps $p_{m}$, we have
\begin{eqnarray}\label{e2.9}
\End' W=\{ f\in \End W\;| \;fp_{m}=0\;\;\;
\mbox{ for } m \mbox{ sufficiently large}\}.
\end{eqnarray}
Clearly, $\End'W$ is a graded subalgebra with respect to 
the grading given in (\ref{eendf1}).
But $\End'W$ does not have a unit.
By using (\ref{e2.9}) it is not hard to see 
that $\End'W$ is an ideal of $gl_{J}(W)$.
(Thus $gl_{J}(W)$ is not a simple algebra.)

We shall need the following simple fact from linear algebra.

\begin{lemma}\label{lfact}
Let $U_{1}, U_{2}$ and $U_{3}$ be finite-dimensional vector spaces 
with $U_{2}\ne 0$.
Then
\begin{eqnarray}\label{eu123}
& &\Hom(U_{1},U_{3})=\Hom (U_{2},U_{3})\cdot \Hom (U_{1},U_{2})\\
& &\hspace{2cm}\left(=\mbox{linear span}
\{fg\;|\;f\in \Hom(U_{2},U_{3}),\;g\in \Hom(U_{1},U_{2})\}\right).
\nonumber
\end{eqnarray}
\end{lemma}

\begin{proof}
 If either $U_{1}=0$ or $U_{3}=0$, it is clear. Suppose $U_{1}\ne 0$
and $U_{3}\ne 0$. 
Let $u_{1},\dots, u_{r}$ be a basis of $U_{1}$ and $v_{1},\dots, v_{s}$
be a basis of $U_{3}$.
For any $1\le i\le r, \;1\le j\le s$, let $f_{ij}$ be the linear
homomorphism
defined by $f_{ij}(u_{t})=\delta_{it}v_{j}$ for $1\le t\le r$.
Let $g\in \Hom(U_{1},U_{2})$ be such that
$g(u_{i})\ne 0$ and $g(u_{k})=0$ for $k\ne i$ and 
let $f\in \Hom (U_{2},U_{3})$ 
be such that $f(g(u_{i}))=v_{j}$.
Then $f_{ij}=fg$. From this (\ref{eu123}) follows immediately.
\end{proof}

\begin{remark}\label{rendf} 
Notice that $p_{m}\in \End'W$ for $m\in {\mathbb Z}$
and that
$$p_{m}fp_{n}\in \Hom(W(n),W(m))\;\left(\subset \End' W\right)$$
for $m,n\in {\mathbb Z}, \;f\in \End'W.$
Then using Lemma \ref{lfact} one can easily show
that $0$ and $\End'W$ are the only ideals of $\End 'W$.
That is, $\End 'W$ is simple.
\end{remark}

In the following we review
a certain formal completion of a ${\mathbb Z}$-graded associative algebra
(cf. [FZ], [KL], [M]).

Let $A=\coprod_{n\in {\mathbb Z}}A_{n}$ be a ${\mathbb Z}$-graded
associative algebra {\em with or without a unit}. For $m,k\in {\mathbb Z}$,
 we set
\begin{eqnarray}
A_{m,k}=\sum_{n>k}A_{m+n}A_{-n}\;\;\;\left(\subset A_{m}\right).
\end{eqnarray}
Then
\begin{eqnarray}\label{ebasic1}
A_{m,k+1}\subset A_{m,k},\;\;\;A_{n}A_{m,k}\subset A_{m+n,k},
\;\;\; A_{m,k}A_{n}\subset A_{m+n,k-n}
\end{eqnarray}
for $m,n,k\in {\mathbb Z}$. Using $A_{m,k}$ for $k\in {\mathbb Z}$ as basic open
neighborhoods of $0\in A_{m}$, we endow $A_{m}$ a topology
with which $A_{m}$ becomes a topological vector space.
Let $\tilde{A}_{m}$ be the completion with respect to this topology
and then set $\tilde{A}=\coprod_{m\in {\mathbb Z}}\tilde{A}_{m}$.
It follows from (\ref{ebasic1}) that the multiplication of 
$A$ is continuous so that $\tilde{A}$ is a ${\mathbb Z}$-graded
topological algebra.

We can explicitly define each $\tilde{A}_{m}$ by using Cauchy sequences
as follows:
An element  $f$ of $\Map({\mathbb N},A_{m})$ is said to be {\em Cauchy} if
for any $k\in {\mathbb Z}$, there exists $r\ge 0$ such that
$f_{i}-f_{j}\in A_{m,k}$ whenever $i,j\ge r$. Then all Cauchy sequences 
form a vector subspace $C(A_{m})$ of $\Map({\mathbb N},A_{m})$.
We define a relation ``$\sim$'' on $C(A_{m})$ such that
$f\sim g$ if and only if for any $k\in {\mathbb Z}$, 
there exists $r\ge 0$ such that
$f_{i}-g_{i}\in A_{m,k}$ whenever $i\ge r$. 
Clearly, ``$\sim$'' is an equivalent relation.
Then define $\tilde{A_{m}}=C(A_{m})/\sim$ and
\begin{eqnarray}\label{e2.12}
\tilde{A}=\coprod_{m\in {\mathbb Z}}\tilde{A_{m}}.
\end{eqnarray}
For $f\in C(A_{m}), \;g\in C(A_{n})$, we define
$fg\in \Map({\mathbb N},A_{m+n})$ by $(fg)_{i}=f_{i}g_{i}$ for $i\in {\mathbb N}$.
It follows from the identity
$$f_{i}g_{i}-f_{j}g_{j}=(f_{i}-f_{j})g_{i}+f_{j}(g_{i}-g_{j})$$
and the property (\ref{ebasic1}) that $fg\in C(A_{m+n})$.
Similarly, if $f, f'\in C(A_{m}), \;g, g'\in C(A_{n})$ and
$f\sim f', g\sim g'$, using the identity 
$$f_{i}g_{i}-f'_{i}g'_{i}=(f_{i}-f_{i}')g_{i}+f_{i}'(g_{i}-g_{i}')$$
and (\ref{ebasic1}) we have $fg\sim f'g'$.
Then we obtain a well defined bilinear map from 
$\tilde{A}_{m}\times \tilde{A}_{n}$ to $\tilde{A}_{m+n}$.
Using linearity we obtain a bilinear multiplication on $\tilde{A}$.
Clearly, this makes $\tilde{A}$ a ${\mathbb Z}$-graded associative algebra
with respect to the grading given in (\ref{e2.12}).

Let $\pi$ be the linear map from $A$ to $\tilde{A}$ such that
$\pi(a)$ is the constant map with value $a$ for $a\in A_{m}, \;m\in {\mathbb Z}$.
Then $\pi$ is a grading-preserving algebra homomorphism. 
For $a\in A_{m},\;m\in {\mathbb Z}$, 
$\pi(a)=0$ if and only if $a\in \cap_{k\in {\mathbb Z}}A_{m,k}$. Therefore
\begin{eqnarray}
\ker \pi=\coprod_{m\in {\mathbb Z}}\left(\cap_{k\in {\mathbb Z}}A_{m,k}\right).
\end{eqnarray}
If $\coprod_{m\in {\mathbb Z}}\left(\cap_{k\in {\mathbb Z}}A_{m,k}\right)=0$, $A$ is 
embedded into $\tilde{A}$ as
a subalgebra through $\pi$.
In general, $\pi$ may not be injective. However, this is not too bad 
when we consider the so-called lower truncated ${\mathbb Z}$-graded 
$A$-modules, which we define next.

\begin{definition}\label{dlower}
{\em A lower truncated ${\mathbb Z}$-graded} $A$-module is an
$A$-module $M$ equipped with a grading
$M=\coprod_{n\in {\mathbb Z}}M(n)$  such that for $m,n\in {\mathbb Z}$,
\begin{eqnarray}
& &A_{m}M(n)\subset M(m+n),\\
& &M(n)=0\;\;\;\mbox{ for }n \mbox{ sufficiently small}.
\end{eqnarray}
Two such $A$-modules are said to be {\em equivalent} if
there is a homogeneous $A$-module isomorphism (of some degree) from one to 
the other. The notion of ${\mathbb N}$-graded $A$-module is 
defined in the obvious way.
\end{definition}

With this definition, by shifting the grading any nonzero lower 
truncated ${\mathbb Z}$-graded $A$-module $M$
is equivalent to an ${\mathbb N}$-graded $A$-module $M=\coprod_{n\ge 0}M(n)$
such that $M(0)\ne 0$.

Let $M=\coprod_{n\in {\mathbb N}}M(n)$ be an ${\mathbb N}$-graded
$A$-module. Then for any $m,n\in {\mathbb Z}$, $A_{m,k}M(n)=0$ for $k$ 
sufficiently large, hence $(\cap_{k\in {\mathbb Z}}A_{m,k})M(n)=0$.
Consequently, $(\ker \pi) M(n)=0$ for each $n$, hence
$(\ker \pi) M=0$. That is,  any ${\mathbb N}$-graded
$A$-module is a natural $A/(\ker \pi)$-module. It follows that
the category of ${\mathbb N}$-graded $A$-modules is 
equivalent to the category of ${\mathbb N}$-graded 
$A/(\ker \pi)$-modules.

We naturally extend the action of $A$ 
on an ${\mathbb N}$-graded $A$-module $M$ to
an action of $\tilde{A}$ as follows:
let $f\in \tilde{A}, \;u\in M$. Since $M$ is ${\mathbb N}$-graded,
$A_{m,k}u=0$ for $k$ sufficiently large, so that
$(f_{i+1}-f_{i})(u)=0$ for $i$ sufficiently large. Then
$f_{r}(u)=f_{r+1}(u)=f_{r+2}(u)=\cdots$ for some $r\ge 0$.
Now we define 
\begin{eqnarray}
f(u)=f_{r}(u)=\lim_{i\rightarrow \infty}f_{i}(u).
\end{eqnarray}
It is routine to check that  $M$ is an $\tilde{A}$-module
with the defined action.
Conversely, any lower truncated ${\mathbb Z}$-graded $\tilde{A}$-module
is a natural lower truncated ${\mathbb Z}$-graded $A$-module through 
the algebra homomorphism $\pi$.

We define the following notion of continuous module (cf. [MP]).

\begin{definition}\label{dcontinuous}
Let $B=\coprod_{m\in {\mathbb Z}}B_{m}$ be a ${\mathbb Z}$-graded
topological associative algebra. A lower truncated ${\mathbb Z}$-graded 
$B$-module $M$ is said to be {\em continuous} if when endowed with
the discrete topology $M$ is a continuous $B$-module in the usual 
sense,  i.e., the action map from $B\times M$ to $M$ 
is continuous.
\end{definition}

\begin{remark}\label{rcontin} 
With this notion, any lower truncated ${\mathbb Z}$-graded $A$-module 
$M$ is
a continuous $A$-module where $A$ is endowed with
the topology defined before.
Indeed, for any $m_{1},m_{2}\in M$, let $a\in A_{m}$ be such that
$am_{1}=m_{2}$ and $A_{m,k}m_{1}=0$ for some $k$. Then
$(a+A_{m,k})m_{1}=m_{2}$. Using this one can easily show that
$M$ is a continuous module.
Furthermore, a lower truncated ${\mathbb Z}$-graded $\tilde{A}$-module, 
or a representation $\rho$ is 
continuous if and only if $\rho$ is the natural extension of
the representation $\rho\pi$ of $A$.
However, {\em the left regular module $B$ may not be continuous
under this definition.} In fact, one can show that $B$ is a 
continuous $B$-module if and only if the topology on 
$B$ as a topological algebra is discrete.
\end{remark}

Now we show that $\tilde{A}=gl_{J}(W)$ with $A=\End'W$.

\begin{proposition}\label{pgl{J}(W)}
Let $W=\coprod_{n\ge 0}W(n)$ be an ${\mathbb N}$-graded vector space
with finite-dimensional homogeneous subspaces. Then
$gl_{J}(W)$ as a ${\mathbb Z}$-graded associative algebra is isomorphic
to $\tilde{A}$ where $A=\End'W$.
\end{proposition}

\begin{proof} We first review a well known fact.
Let $U=\coprod_{n\in {\mathbb Z}}U(n)$ be a ${\mathbb Z}$-graded vector space.
Endow $U$ with a topology by using $u+\sum_{n\ge k}U(n)$ for 
$u\in U, \;k\in {\mathbb Z}$ as a base of open sets in $U$. $U$ is a
Hausdorff topological vector space because 
$\cap_{k\in {\mathbb Z}}U(k)=0$. It is well known that
$\prod_{n\in {\mathbb Z}}U(n)$ is the completion of the topological 
vector space $U$ defined above. 
In view of (\ref{e2.5}), for each $m\in {\mathbb Z}$,
 $gl_{J}(W)_{m}$ is the completion
of $A_{m}$ with the topology
on $A_{m}$ defined by using
$\coprod_{n>k}\Hom(W(n),W(m+n))$ $(k\ge 0)$ as basic
open neighborhoods of $0$.
Now we show that this topology on $A_{m}$ is the same
as the one defined by using $A_{m,k}$ for $k\ge 0$
as basic neighborhoods of $0$.
If $W=0$, it is clear. Now we
assume $W=\coprod_{n\ge 0}W(n)$ with $W(0)\ne 0$.
By Lemma \ref{lfact} we have
$$\Hom(W(n),W(m+n))=\Hom(W(0),W(m+n))\Hom(W(n),W(0))
\subset A_{m+n}A_{-n}.$$
On the other hand, 
\begin{eqnarray*}
A_{m+n}A_{-n}&=&\sum_{r,s}\Hom(W(r),W(m+n+r))\Hom(W(s),W(s-n))\nonumber\\
&\subset&\sum_{s\ge n}\Hom(W(s),W(m+s))
\end{eqnarray*}
because $W(s-n)=0$ for $s<n$.
Thus
\begin{eqnarray}
A_{m,k}=\coprod_{s>k}\Hom(W(s),W(m+s)).
\end{eqnarray}
Therefore, $gl_{J}(W)$ is the completion of $A=\End'(W)$ 
with the given topology. Furthermore, an element $\sum_{n}a_{n}$
of $gl_{J}(W)_{m}$ with $a_{n}\in \Hom(W(n),W(m+n))$ 
corresponds to the Cauchy sequence $\{\alpha_{n}\}$ with
$\alpha_{n}=\sum_{i=0}^{n}a_{i}$. From this, the ${\mathbb Z}$-graded 
associative algebra $gl_{J}(W)$ is the same as the ${\mathbb Z}$-graded 
associative algebra $\tilde{A}$.
\end{proof}

\begin{definition}\label{dtopgl{J}(W)}
We define $gl_{J}(W)$ to be the topological algebra
with the topology obtained by identifying $gl_{J}(W)$ with $\tilde{A}$
where $A=\End'W$.
\end{definition}

In view of this, a lower truncated ${\mathbb Z}$-graded continuous 
$gl_{J}(W)$-module
amounts to a lower truncated ${\mathbb Z}$-graded $\End'W$-module.
Clearly, $W$ is a continuous $gl_{J}(W)$-module.

Define an element $d$ of $gl_{J}(W)_{0}$ by
\begin{eqnarray}
d|_{W(n)}=n\cdot {\rm id}_{W(n)}\;\;\;\;\mbox{ for }n\in {\mathbb Z}.
\end{eqnarray}
Then
\begin{eqnarray}
da-ad=na\;\;\;\;\mbox{ for }a\in gl_{J}(W)_{n}, \;n\in {\mathbb Z}.
\end{eqnarray}
We have the following result.

\begin{proposition}\label{pamodule}
The natural module $W$ is a continuous irreducible $gl_{J}(W)$-module.
\end{proposition}

\begin{proof} With the element $d$ of $gl_{J}(W)$,
it follows that any $gl_{J}(W)$-submodule of $W$
is automatically graded. Since 
for any $k\in {\mathbb N}$, $\oplus_{n\le k}W(n)$ is an irreducible
$\End (\oplus_{n\le k}W(n))$-module and
$$\End (\oplus_{n\le k}W(n))\subset \End' W\subset gl_{J}(W),$$
it follows that $W$ is an irreducible ${\mathbb N}$-graded 
$gl_{J}(W)$-module.
\end{proof}

\section{Rationality of $gl_{J}(W)$}
As before, $W$ will be an ${\mathbb N}$-graded space with finite-dimensional 
homogeneous subspaces.
Our goal of this section is to prove the following:

\begin{theorem}\label{t1}
The natural module $W$ is the unique 
continuous irreducible ${\mathbb N}$-graded
$gl_{J}(W)$-module up to equivalence and any
continuous ${\mathbb N}$-graded $gl_{J}(W)$-module is a direct 
sum of some copies of $W$.
\end{theorem}

We shall prove this theorem as an application of
a slightly more general result. 
The assertions of Theorem \ref{t1} are analogues
of those for a (finite-dimensional) matrix algebra.
Motivated by 
the notion of rationality of vertex operator algebras
(cf. [Z], [DLM2]) we define the following notion:

\begin{definition}\label{drationalA} 
A ${\mathbb Z}$-graded topological associative 
algebra $B$ is said to be {\em rational} if
there are only finitely many irreducible 
${\mathbb N}$-graded continuous $B$-modules up to equivalence and
any ${\mathbb N}$-graded continuous
$B$-module is a direct sum of irreducible graded modules
with finite-dimensional homogeneous subspaces.
\end{definition}

Any (finite-dimensional) semisimple algebra $A$ (over ${\mathbb C}$)
are rational where $A=A_{0}$. Thus the notion of rationality
is a generalization of the classical notion of semisimplicity.
Theorem \ref{t1} implies that $gl_{J}(W)$ is rational.
It is clear that the direct sum of (finitely many) rational 
algebras with the product topology 
are still rational. Then we immediately have:

\begin{corollary}\label{cdirectsum}
Let $W_{1},\dots, W_{r}$ be ${\mathbb N}$-graded
vector spaces with
homogeneous subspaces being  finite-dimensional. Then
$gl_{J}(W_{1})\oplus \cdots \oplus gl_{J}(W_{r})$ is a
rational ${\mathbb Z}$-graded topological algebra.
\end{corollary}

To prove Theorem \ref{t1} we shall use
an analogue of Zhu's $A(V)$-theory (cf. [Z], [DLM2-3]) from
vertex operator algebra theory.

Now we consider a ${\mathbb Z}$-graded topological associative algebra  
$B=\coprod_{m\in {\mathbb Z}}B_{m}$ with unit. Recall that for $k\ge 0$,
$$B_{0,k}=\sum_{n> k}B_{n}B_{-n}\;\;\left(\subset B_{0}\right).$$
Clearly, each $B_{0,k}$ is a two-sided ideal of $B_{0}$ and
so is the closure $\overline{B_{0,k}}$.
Set
\begin{eqnarray}
T_{k}(B)=B_{0}/\overline{B_{0,k}}.
\end{eqnarray}
Then any $B_{0}$-module $U$ with $\overline{B_{0,k}}U=0$ is a natural 
$T_{k}(B)$-module
and on the other hand, any $T_{k}(B)$-module is a natural $B_{0}$-module.
If $M=\coprod_{n\ge 0}M(n)$ is an ${\mathbb N}$-graded continuous $B$-module,
then $M(n)$ is a natural $T_{k}(B)$-module for $k\ge n$ because 
$\overline{B_{0,k}}M(n)=0$.

Set
\begin{eqnarray}
B_{\le 0}=\coprod_{n\le 0}B_{n},
\end{eqnarray}
a graded subalgebra of $B$.
Let $U$ be a $T_{0}(B)$-module. We consider $U$ as a $B_{\le 0}$-module
with $B_{n}U=0$ for $n<0$. Form the induced $B$-module
\begin{eqnarray}
M(U)=B\otimes_{B_{\le 0}}U.
\end{eqnarray}
Then $M(U)$ is an ${\mathbb N}$-graded $B$-module with
$M(U)(0)=U$.
Let $J(U)$ be the sum of all graded submodule
$N$ of $M(U)$ with $N(0)=0$ and set
\begin{eqnarray}
L(U)=M(U)/J(U).
\end{eqnarray}
Then $L(U)$ is an ${\mathbb N}$-graded $B$-module with $L(U)(0)=U$
such that $L(U)$ as a $B$-module is generated by $U$ and that
for any nonzero graded submodule $N$ of $L(U)$, $N(0)\ne 0$.
Consequently, if $U$ is an irreducible $T_{0}(B)$-module, $L(U)$ is
an irreducible ${\mathbb N}$-graded $B$-module. (Here we do not claim that
$L(U)$ is a continuous $B$-module.)

\begin{lemma}\label{ltkb}
Let $W_{i}=\coprod_{n\in {\mathbb N}}W_{i}(n)$ $(i=1,2)$ be ${\mathbb N}$-graded
irreducible continuous $B$-modules with $W_{i}(0)\ne 0$. Then
$W_{1}\cong W_{2}$ as a $B$-module if and only if $W_{1}(0)\cong W_{2}(0)$
as a $T_{0}(B)$-module.
\end{lemma}

\begin{proof} We only need to prove that $W_{1}\cong W_{2}$ 
if $W_{1}(0)\cong W_{2}(0)$ as a $T_{0}(B)$-module.
By the universal property of induced modules,
$W_{1}$ and $W_{2}$ are natural quotient $B$-modules of
$M(W_{1}(0))$ and $M(W_{2}(0))$ by $J(W_{1}(0))$ and $J(W_{2}(0))$,
respectively.
Since $W_{1}(0)$ and $W_{2}(0)$ are equivalent $B_{\le 0}$-modules,
$M(W_{1}(0))$ and $M(W_{2}(0))$ are equivalent ${\mathbb N}$-graded
$B$-modules. 
Then it follows that $W_{1}$ and $W_{2}$
are equivalent $B$-modules.
\end{proof}

The following is a sufficient condition for $B$ to be rational:

\begin{proposition}\label{pAfact}
Let $B=\coprod_{n\in {\mathbb Z}}B_{n}$ be a ${\mathbb Z}$-graded 
topological
associative algebra equipped with an element $d$ of $B_{0}$ such
that
\begin{eqnarray}
db-bd=nb\;\;\;\;\mbox{ for }b\in B_{n}, \;n\in {\mathbb Z}.
\end{eqnarray}
Then $B$ is rational if the following conditions
are satisfied:

(1) for every $k\ge 0$, $T_{k}(B)$ is (finite-dimensional) semisimple;
 
(2) for any two inequivalent
 irreducible $T_{0}(B)$-modules $U_{1}$ and $U_{2}$,
$d_{1}-d_{2}\notin {\mathbb Z}$, where $d$ acts on $U_{i}$ 
as a scalar $d_{i}$.

(3) for any ${\mathbb N}$-graded continuous $B$-module
$W=\coprod_{n\in {\mathbb N}}W(n)$,
$(B_{-n}B_{n})W(0)\ne 0$ if $B_{n}W(0)\ne 0$ for $n\ge 1$.
\end{proposition}

\begin{proof} Let us assume (1), (2) and (3). 

Claim 1: Let  $M=\coprod_{n\in {\mathbb N}}M(n)$ be
any ${\mathbb N}$-graded continuous $B$-module such that
$M(0)$ is an irreducible $B_{0}$-module and $M=BM(0)$. Then
$M$ is an irreducible $B$-module.

Since $T_{0}(B)$ is finite-dimensional and $M(0)$ is an
irreducible $T_{0}(B)$-module, $M(0)$ is finite-dimensional.
Then being a central element 
of $B_{0}$, $d$ acts as a scalar $\lambda$ 
on $M(0)$. From $BM(0)=M$, we get
$M(n)=B_{n}M(0)$ for $n\ge 0$. Then 
$$dw=(\lambda+n)w\;\;\;\mbox{ for }w\in M(n), \;n \ge 0.$$
(We call $\lambda$ the {\em lowest weight} of $M$.)
Consequently, any submodule of $M$ is graded.
Let $N=\oplus_{n\in {\mathbb N}}N(n)$ 
be any nonzero submodule of $M$. We shall prove that
$N(0)=M(0)$, which implies that $T=M$ because $M=BM(0)$.
Let $k$ be a nonnegative integer such that 
$N(k)\ne 0$.
Since $B_{k}M(0)=M(k)\supset N(k)\ne 0$, by hypothesis (3)
 $B_{-k}B_{k}M(0)=M(0)$, hence
\begin{eqnarray}\label{e*}
B_{k}B_{-k}M(k)=B_{k}B_{-k}B_{k}M(0)=B_{k}M(0)=M(k)\ne 0.
\end{eqnarray}
Since $M(k)$ is a completely reducible $B_{0}$-module, from (\ref{e*})
 there is
an irreducible $T_{k}(B)$-submodule $S$ of $M(k)$ such that
$B_{k}B_{-k}S\ne 0$.
Then $B_{-k}S=M(0)$ because $B_{-k}S$ is a nonzero $B_{0}$-submodule of 
$M(0)$ 
and $M(0)$ is irreducible. Consequently, $B_{k}B_{-k}S=B_{k}M(0)=M(k)$.
On the other hand, $(B_{k}B_{-k})S\subset S$. Thus $M(k)=S$
is an irreducible $T_{k}(B)$-module, hence $N(k)=M(k)$.
If $k=0$, we have $N(0)=M(0)$. If $k>0$, we have 
$B_{-k}M(k)=B_{-k}B_{k}M(0)=M(0)$, 
so that
$$N(0)\supset B_{-k}N(k)=B_{-k}M(k)=M(0).$$
Then $N(0)=M(0)$. Thus $N=M$.
Therefore $M$ is irreducible.

Claim 2: any ${\mathbb N}$-graded continuous $B$-module
$W=\coprod_{n\in {\mathbb N}}W(n)$ is 
a direct sum of irreducible graded $B$-modules.

Let $W^{o}$ be the sum of
all irreducible graded submodules of $W$.
We must prove that $W=W^{o}$.
Since $T_{0}(B)$ is semisimple, 
$W(0)$ as a $B_{0}$-module is completely reducible.
By Claim 1 the $B$-submodule generated by $W(0)$ is completely reducible,
so that $W(0)\subset W^{o}(0)$.
Suppose that for some nonnegative integer $k$, 
$W(n)=W^{o}(n)$ for $n\le k$.
Write $W(k+1)=W^{o}(k+1)\oplus P$, where $P$ is a $B_{0}$-submodule of
$W(k+1)$. Then
$$B_{-n}P\subset W(k+1-n)=W^{o}(k+1-n)\;\;\;\mbox{ for }n\ge 1,$$
hence
$$B_{n}B_{-n}P\subset B_{n}W^{o}(k+1-n)\subset W^{o}(k+1)
\;\;\;\mbox{ for }n\ge 1.$$
On the other hand, $B_{n}B_{-n}P\subset B_{0}P=P$.
Thus $B_{n}B_{-n}P=0$ for $n\ge 1$, hence $B_{0,0}P=0$.
Since $W$ is continuous, $\overline{B_{0,0}}P=0$, so that
$P$ is a $T_{0}(B)$-module. Now we claim that $P=0$.
If $P\ne 0$, then let $U$ be an irreducible ($T_{0}(B)$-) 
$B_{0}$-submodule of $P$. Then $d=\lambda$ on $U$ for some
 $\lambda\in {\mathbb C}$.
If $B_{-n}U=0$ for all $n\ge 1$, 
by Claim 1, $U$ generates an irreducible
graded $B$-submodule of $W$. By the definition of $W^{o}$, 
$U\subset W^{o}$. It is a contradiction. 
If $B_{-n}U\ne 0$ for some $n\ge 1$, then 
$\lambda'-\lambda+n\in {\mathbb N}$, where $\lambda'$ is the
lowest weight of some irreducible ${\mathbb N}$-graded continuous
$B$-module because
$$B_{-n}U\subset W(k-n)=W^{o}(k-n)\subset W^{o}.$$
This contradicts Hypothesis (2). Thus $P=0$,
hence $W(k+1)=W^{o}(k+1)$.
Therefore $W=W^{o}$ by induction.
\end{proof}


To apply Proposition \ref{pAfact} for $B=gl_{J}(W)$ we need to calculate
$T_{k}(B)$ for $k\ge 0$.
Recall that for $k\ge 0$,
$T_{k}(gl_{J}(W))=gl_{J}(W)_{0}/\overline{gl_{J}(W)_{0,k}}$, where
\begin{eqnarray}
\overline{gl_{J}(W)_{0,k}}=\prod_{n\ge k}gl_{J}(W)_{m+n}gl_{J}(W)_{-n}
\;\;\left(\subset gl_{J}(W)_{0}\right).
\end{eqnarray}

\begin{lemma}\label{lawk}
For $k\ge 0$, we have
\begin{eqnarray}\label{e2.27}
T_{k}(gl_{J}(W))=\End W(0)\oplus \cdots \oplus \End W(k).
\end{eqnarray}
In particular, $T_{k}(gl_{J}(W))$ is semisimple.
\end{lemma}

\begin{proof} First, we have
\begin{eqnarray}
gl_{J}(W)_{0}=\prod_{n\ge 0}\End W(n).
\end{eqnarray}
By Lemma \ref{lfact}, for $n\in {\mathbb N}$,
\begin{eqnarray}\label{e2.33}
& &\End W(n)=\Hom (W(0),W(n))\Hom (W(n),W(0))\\
& &\subset (\End'W)_{n}(\End'W)_{-n}.\nonumber
\end{eqnarray}
Let $\pi_{k}$ be the projection of $gl_{J}(W)_{0}$ onto
$\oplus_{n=0}^{k}\End W(n)$. (This is an algebra homomorphism.)
Now we calculate the kernel of $\pi_{k}$. By definition, $\pi_{k}(f)=0$
if and only $f_{0}=f_{1}=\cdots =f_{k}=0$, which is equivalent
to that $f\in \prod_{n> k}\End W(n)$.
Then it follows from (\ref{e2.33}) that
$\ker \pi_{k}\subset \overline{gl_{J}(W)_{0,k}}$.
Clearly, $\overline{gl_{J}(W)_{0,k}}\subset \ker \pi_{k}$.
Thus, $\ker \pi_{k}=\overline{gl_{J}(W)_{0,k}}$.
Therefore, $\pi_{k}$ gives rise to an algebra isomorphism
for (\ref{e2.27}).
\end{proof}

{\bf Proof of Theorem \ref{t1}:}
It follows from Lemmas \ref{lawk} and \ref{ltkb} that
$W$ is the unique continuous irreducible ${\mathbb N}$-graded 
$gl_{J}(W)$-module up to equivalence.
In view of Proposition \ref{pAfact}
and Lemma \ref{lawk} it suffices to prove that
$gl_{J}(W)_{-n}gl_{J}(W)_{n}W(0)\ne 0$ if $W(n)\ne 0$.
{}From Lemma \ref{lfact},
\begin{eqnarray*}
& &gl_{J}(W)_{-n}gl_{J}(W)_{n}|_{W(0)}\\
& &\supset \Hom(W(n),W(0))\Hom(W(0),W(n))|_{W(0)}
=\End W(0)\ne 0.
\end{eqnarray*}
Then Proposition \ref{pAfact} applies.

\section{The analogue of the Burnside theorem}

Throughout this section, $A$ will be a ${\mathbb Z}$-graded
associative algebra and $\rho$ will be an irreducible ${\mathbb N}$-graded 
representation of $A$ on $W=\oplus_{n\ge 0}W(n)$ where $W(0)\ne 0$ and
all $W(n)$ are finite-dimensional.

Here we {\em do not assume} that $A$ has a unit. By an irreducible
$A$-module $M$ we mean that $0$ and $M$ are the only submodules
and $AM\ne 0$.

\begin{lemma}\label{l3.1} We have: (a)
$W(n)$ is an irreducible $A_{0}$-module if $W(n)\ne 0$;
(b) $W(r)$ and $W(s)$ are inequivalent $A_{0}$-modules
if $r\ne s, W(r)\ne 0$ and $W(s)\ne 0$.
\end{lemma}

\begin{proof} Set
\begin{eqnarray}
ann_{W}(A)=\{w\in W| Aw=0\}.
\end{eqnarray}
Then $ann_{W}(A)$ is an $A$-submodule of $W$. Since $AW\ne 0$,
$ann_{W}(A)\ne W$. By the irreducibility of $W$, we have $ann_{W}(A)=0$,
hence $Aw\ne 0$ for $0\ne w\in W$. By the irreducibility of $W$ again,
 we have
$Aw=W$ for any $0\ne w\in W$. Thus $A_{m}w=W(m+n)$ for any 
$m\in {\mathbb Z}$ if $0\ne w\in W(n)$. Consequently, $W(n)$ is an 
irreducible $A_{0}$-module and $A_{m}W(n)=W(m+n)$ for any 
$m\in {\mathbb Z}$ if $W(n)\ne 0$. If $W(r)\ne 0, W(s)\ne 0$ for
some $r>s$, then 
$$(A_{r}A_{-r})W(r)=A_{r}W(0)=W(r)\ne 0,\;\;\;\; 
(A_{r}A_{-r})W(s)=A_{r}(A_{-r}W(s))=0.$$
 Therefore,
$W(r)$ and $W(s)$ are inequivalent $A_{0}$-modules.
\end{proof}

\begin{remark}\label{rclass} 
Here we mention that the standard Burnside theorem
still holds for algebras without a unit.
Let $A$ be an algebra without a unit and $\rho$ be a
finite-dimensional irreducible representation of $A$ on $U$.
First, extend $A$ to be an algebra
$\bar{A}=A\oplus {\mathbb C}1$ with a unit (cf. [J1]). Then 
$U$ is an irreducible $\bar{A}$-module.
By the standard Burnside theorem,
$\rho(A)+{\mathbb C}1=\rho(\bar{A})=\End U$.
This implies that $\rho(A)$ is an ideal of $\End U$. Since
$\End U$ is simple and $\rho(A)\ne 0$, we have $\rho(A)=\End U$.
Furthermore, the following Chinese Remainder theorem also holds:
if $U=U_{1}\oplus \cdots \oplus U_{r}$ is a faithful
$A$-module where $U_{i}$ are inequivalent finite-dimensional irreducible 
$A$-modules, then $A=\End U_{1}\oplus \cdots \oplus \End U_{r}$.
\end{remark}

Since each $W(n)$ is an $A_{0}$-module, the decomposition 
$W=\coprod_{n\ge 0}W(n)$ gives rise to
an algebra homomorphism $\rho_{k}$ from
$A_{0}$ to $\End W(k)$ for each $k\ge 0$.
For $k\ge 0$, set
\begin{eqnarray}
W^{(k)}=\oplus_{n=0}^{k}W(n)
\end{eqnarray}
and $\rho^{(k)}=\rho_{0}+\cdots +\rho_{k}$, an algebra 
homomorphism from $A_{0}$ to $\End W^{(k)}$.

\begin{lemma}\label{ltka} For $k\ge 0$, we have
\begin{eqnarray}
& &\rho^{(k)}(A_{0})=\End W(0)\oplus \cdots \oplus \End W(k),
\label{etka1}\\
& &\rho^{(k)}(A_{k}A_{-k})=\End W(k).\label{etka2}
\end{eqnarray}
\end{lemma}

\begin{proof} (\ref{etka1}) directly follows from Lemma \ref{l3.1} 
and Remark \ref{rclass}.
Since $A_{k}A_{-k}$ is an ideal of $A_{0}$, it follows from
(\ref{etka1}) that $\rho^{(k)}(A_{k}A_{-k})$ is 
an ideal of $\oplus_{j=0}^{k}\End W(j)$. 
Since $A_{k}A_{-k}W(n)=0$ for $n<k$, $\rho^{(k)}(A_{k}A_{-k})$ is an 
ideal of $\End W(k)$.
If $W(k)=0$, there is nothing to prove.
Suppose $W(k)\ne 0$. Since $A_{k}A_{-k}W(k)=A(k)W(0)=W(k)\ne 0$,
$\rho^{(k)}(A_{k}A_{-k})\ne 0$, hence 
$\rho^{(k)}(A_{k}A_{-k})=\End W(k)$ because $\End W(k)$ is simple.
\end{proof}

We have the following analogue of the Burnside theorem (cf. [M]):

\begin{theorem}\label{tAW} Let $A=\coprod_{m\in {\mathbb Z}}A_{m}$ be 
a ${\mathbb Z}$-graded associative algebra
and let $W$ be an irreducible ${\mathbb N}$-graded $A$-module with 
$W(0)\ne 0$
and with $W(n)$ being finite-dimensional for each $n$.
Then $\tilde{\rho}(\tilde{A})=gl_{J}(W)$.
\end{theorem}

\begin{proof} Let $\psi\in gl_{J}(W)_{n}$.
As usual, we shall use $a$ for $\rho(a)$.
Recall that
$gl_{J}(W)_{n}=\prod_{k\ge 0}\Hom(W(k),W(n+k))$.
Let $\psi_{k}$ be the 
projection of $\psi$ in the subspace $\Hom (W(k),W(n+k))$ 
for $k\in {\mathbb N}$.
Now we are going to find a sequence of elements 
$a_{0}, a_{1},\dots$ of $A$ such that
\begin{eqnarray}
& &a_{r}\in A_{n+r}A_{-r},\label{en=+1}\\
& &\psi_{0}+\cdots +\psi_{r}=(a_{0}+\cdots +a_{r})|_{W^{(r)}}
\label{en=+2}
\end{eqnarray}
for $r\ge 0$. Then we will have
$a_{0}+a_{1}+\cdots \in \tilde{A}$ and
\begin{eqnarray}
\psi=\tilde{\rho}(a_{0}+a_{1}+\cdots).
\end{eqnarray}

Fact: for any
 $\phi\in \Hom(W(0),W(m))$ $(m\ge 0)$,
there exists $a(\phi)\in A_{m}A_{0}$ such that $a(\phi)|_{W(0)}=\phi$.

Since $\rho^{(m)}(A_{m}A_{-m})=\End W(m)$ (by Lemma \ref{ltka}),
there are elements
$b_{m}^{(i)}\in A_{m}, \;c_{m}^{(i)}\in A_{-m}$ such that
\begin{eqnarray}
1|_{W(m)}=\sum_{i}b_{m}^{(i)}c_{m}^{(i)}|_{W(m)}.
\end{eqnarray}
There exist $d_{i}\in A_{0}$ such that
$c_{m}^{(i)}\phi=d_{i}|_{W(0)}$ because
$c_{m}^{(i)}\phi\in \End W(0)=\rho_{0}(A_{0})$ (by Lemma \ref{ltka}).
Now set
\begin{eqnarray}
a(\phi)=\sum_{i}b_{m}^{(i)}d_{i}\in A_{m}A_{0}.
\end{eqnarray}
Then $a(\phi)|_{W(0)}=\sum_{i}b_{m}^{(i)}c_{m}^{(i)}\phi=\phi$. 

If $n\ge 0$, by taking $m=n,\; \phi=\psi_{0}$,
then setting $a_{0}=a(\phi)\in A_{n}A_{0}$, we have
$$a_{0}|_{W(0)}=\psi_{0}.$$
If $n<0$, we have $W(n)=W(n+1)=\cdots =W(-1)=0$, hence
$\psi_{0}=\cdots =\psi_{-n-1}=0$. In this case,
we simply take $a_{0}=\cdots =a_{-n-1}=0$.

Suppose that we have already found $a_{0},\cdots, a_{k}$ such that
(\ref{en=+1}) and (\ref{en=+2}) hold for $0\le r\le k$.
Since
$$(\psi_{k+1}-a_{0}-\cdots -a_{k})b_{k+1}^{(i)}|_{W(0)}
\in \Hom(W(0),W(n+k+1)),$$
by the previous fact, there are
$e_{i}\in A_{n+k+1}A_{0}$ such that
$$(\psi_{k+1}-a_{0}-\cdots -a_{k})b_{k+1}^{(i)}|_{W(0)}=e_{i}|_{W(0)}.$$
Set
\begin{eqnarray}\label{eak+1}
a_{k+1}=\sum_{i}e_{i}c_{k+1}^{(i)}\;\;\in A_{n+k+1}A_{-k-1}.
\end{eqnarray}
Then
\begin{eqnarray}
a_{k+1}|_{W(k+1)}
&=&\sum_{i}e_{i}c_{k+1}^{(i)}|_{W(k+1)}\nonumber\\
&=&\sum_{i}e_{i}|_{W(0)}c_{k+1}^{(i)}|_{W(k+1)}\nonumber\\
&=&\sum_{i}(\psi_{k+1}-a_{0}-\cdots -a_{k})
b_{k+1}^{(i)}c_{k+1}^{(i)}|_{W(k+1)}\nonumber\\
&=&\psi_{k+1}-(a_{0}+\cdots +a_{k})|_{W(k+1)}.
\end{eqnarray}
That is,
\begin{eqnarray}
\psi_{k+1}=(a_{0}+\cdots +a_{k}+a_{k+1})|_{W(k+1)}.
\end{eqnarray}
Since $A_{-k-1}W^{(k)}=0$, 
$a_{k+1}W^{(k)}=0$ (by (\ref{eak+1})).
Thus
$$\psi_{0}+\cdots +\psi_{k}+\psi_{k+1}
=(a_{0}+\cdots +a_{k}+a_{k+1})|_{W^{(k+1)}}.$$
In this way we obtain a sequence of elements $a_{0},a_{1},\dots$
with the desired properties.
This completes the proof. 
\end{proof}

For any $A$-module $W$, set
\begin{eqnarray}
ann_{A}(W)=\{ a\in A\;|\; aW=0\}.
\end{eqnarray}
The following is an immediate consequence of Theorem \ref{tAW} 
and the Chinese Remainder Theorem:

\begin{corollary}\label{chom}
Let $W=W_{1}\oplus \cdots \oplus W_{r}$ be a (semisimple)
${\mathbb N}$-graded $A$-module
with $\rho$ being the representation homomorphism from $A$ to
$gl_{J}(W)$, where $W_{i}$ are inequivalent irreducible $A$-modules.
Then $\rho$ gives rise to an algebra homomorphism $\tilde{\rho}$ from 
$\tilde{A}$ onto
$\prod_{i=1}^{r}gl_{J}(W_{i})$ where
\begin{eqnarray}
\ker \tilde{\rho}=\cap_{i=1}^{r}ann_{\tilde{A}}(W_{i}).
\end{eqnarray}
\end{corollary}

\begin{remark}\label{rma} 
The assertion of Theorem \ref{tAW} is different from
saying that $\rho(A)$ is a dense subspace of $gl_{J}(W)$.
The last 
statement can be proved by using the classical density theorem as follows:
for any $f\in gl_{J}(W)_{n}$, by the classical density theorem [J2], 
for each $k\ge 0$,
there exists an $a_{k}\in A$ such that
$f|_{W^{(k)}}=a_{k}|_{W^{(k)}}$. Since $(f-a_{r})|_{W^{(k)}}=0$
for $r>k$, we have
\begin{eqnarray*}
& &(f-a_{r})|_{W^{(k)}}\in \prod_{m\ge k+1}\Hom(W(m),W(m+n))\\
& &\subset \prod_{m\ge k+1}\Hom(W(0),W(m+n))\Hom(W(m),W(0))\\
& &\subset \overline{\sum_{m\ge k+1}(\End'W)_{m+n}(\End'W)_{-m}}.
\end{eqnarray*}
Thus $\lim_{k\rightarrow \infty}\rho(a_{k})=f$.
\end{remark}

The following is analogous to an exercise in finite-dimensional algebra 
theory [J2]
with which $\End U$ can be alternatively proved to be simple 
(for a finite-dimensional $U$).

\begin{corollary}\label{csimplicity}
Let $I$ be a nonzero ideal of $gl_{J}(W)$ and let $\rho$ be
the representation map of $I$ on $W$. 
Then $\tilde{\rho}(\tilde{I})=gl_{J}(W)$.
\end{corollary}

\begin{proof} If $I$ acts irreducibly on $W$, by Theorem \ref{tAW} 
$\tilde{\rho}(\tilde{I})=gl_{J}(W)$.
Now it suffices to prove that $W$ is an irreducible $I$-module.

Since $I$ is a (left) ideal, for every $w\in W$,
$I\cdot w$ is a $gl_{J}(W)$-submodule
of $W$, hence $I\cdot w=W$ if $I\cdot w\ne 0$. Now 
it suffices to prove that
$I\cdot w\ne 0$ for every nonzero $w\in W$. 
Since $I$ is a (right) ideal, $ann_{W}(I)$ is a $gl_{J}(W)$-submodule 
of $W$, hence $ann_{W}(I)=0$ or $W$. 
We must have $ann_{W}(I)=0$ because $W$ is faithful and $I\ne 0$.
Thus $I\cdot w\ne 0$ if $w\ne 0$. 
\end{proof}


\section{The analogue of the Brauer theorem}
Throughout this section $W=\coprod_{n\in {\mathbb N}}W(n)$ will be 
an ${\mathbb N}$-graded vector space 
such that $W(0)\ne 0$ and all $W(n)$ are finite-dimensional, and
$A$ will be a graded subalgebra 
containing $1$ of $gl_{J}(W)$ such that $W$ is a semisimple $A$-module.
{}From the assumption, Schur's Lemma holds for each of the 
irreducible $A$-submodules 
of $W$. Set
\begin{eqnarray}
C(A)=\{ b\in gl_{J}(W)\;|\; ab=ba \;\;\;\mbox{ for }a\in A\}.
\end{eqnarray}
Then $C(A)$ is a graded subalgebra of $gl_{J}(W)$.

Let $U_{i}$ ($i\in I$) be a complete set of representatives of
equivalence classes of irreducible $A$-submodules of $W$. Then
we have a canonical decomposition
\begin{eqnarray}
W=\coprod_{i\in I}W_{i},
\end{eqnarray}
where $W_{i}$ is the sum of all graded $A$-submodules of $W$
which are equivalent to $U_{i}$ for $i\in I$. 
Furthermore, each $W_{i}$ can be canonically 
decomposed as
\begin{eqnarray}
W_{i}=U_{i}\otimes V_{i},
\end{eqnarray}
where $V_{i}=\Hom_{A}(U_{i},W)$ for $i\in I$.
(Schur's Lemma was used here.) The following proposition
is a simple analogue of the classical duality result (cf. [DLM1], [H]).

\begin{proposition}\label{p2} 
The spaces $V_{i}\;\left(=\Hom_{A}(U_{i},W)\right)$ $(i\in I)$
 are inequivalent 
irreducible graded $C(A)$-modules, so that $W$ is a semisimple
$C(A)$-module.
\end{proposition}

\begin{proof} (1) For any fixed $i\in I$, let
$f,g\in \Hom_{A}(U_{i},W)$ be any two linearly independent 
homogeneous elements of degrees $r$ and $s$, respectively.
We decompose $W$ as
$$W=U_{i}\otimes {\mathbb C}f\oplus U_{i}\otimes {\mathbb C}g\oplus T,$$
where $T$ is an $A$-submodule of $W$.
Define $b\in \End W$ by
\begin{eqnarray}\label{edefb}
b(u\otimes f+v\otimes g+t)=u\otimes g+v\otimes f+t
\end{eqnarray}
for $u,v\in U_{i}, t\in T$. It follows directly that
$ab=ba$ for all $a\in A$. From
(\ref{edefb}) we have $p_{m}bp_{n}=0$ if $|m-n|>|r-s|$. 
Thus $b\in gl_{J}(W)$, hence
$b\in C(A)$. Clearly, $bf=g$. Therefore $\Hom_{A}(U_{i},W)$ is 
an irreducible graded $C(A)$-module.

(2) For $i\ne j\in I$, let
$$W=U_{i}\oplus U_{j}\oplus S,$$
where $S$ is an $A$-submodule of $W$. Define
$e_{i}\in \End W$ by
\begin{eqnarray}
e_{i}|_{U_{i}}=1, \;\;e_{i}|_{U_{j}}=0, \;\;e_{i}|_{S}=0.
\end{eqnarray}
Then $e_{i}\in gl_{J}(W)$ and $e_{i}\in C(A)$. 
Since $e_{i}U_{i}=U_{i}$ and $e_{i}U_{j}=0$, $U_{i}$ and $U_{j}$
are inequivalent $C(A)$-modules. 
\end{proof}

We shall continue studying $C(A)$.



\begin{proposition}\label{pdual1}
Let $A$ be a graded subalgebra containing $1$ of $gl_{J}(W)$ such that
$W$ is a semisimple $A$-module. Then
\begin{eqnarray}
& &C(A)=\coprod_{i\in I}\left({\mathbb C}\otimes gl_{J}(V_{i})\right),
\label{eC(A)1}\\ 
& &C(C(A))=\coprod_{i\in I}\left(gl_{J}(U_{i})\otimes {\mathbb C}\right)
=\rho(\tilde{A}).
\end{eqnarray}
\end{proposition}

\begin{proof} Clearly,
\begin{eqnarray}
C(A)\supset \coprod_{i\in I}{\mathbb C}\otimes gl_{J}(V_{i}).
\end{eqnarray}
Conversely, let $b\in C(A)$. Then $bW_{i}\subset W_{i}$ for $i\in I$,
that is, 
$$b(U_{i}\otimes V_{i})\subset (U_{i}\otimes V_{i}).$$
For each fixed $i$, let $v_{j}$ $(j\in J)$ be a basis for $V_{i}$
and let $v$ be any nonzero vector in $V_{i}$.
Define linear endomorphisms $\psi_{j}$ of $U_{i}$ by
\begin{eqnarray}
b(u\otimes v)=\sum_{j\in J}\psi_{j}(u)\otimes v_{j}
\end{eqnarray}
for $u\in U_{i}$. Since $b\in C(A)$, each $\psi_{j}$ is an $A$-endomorphism.
By Schur's Lemma we have $\psi_{j}=\alpha_{j}\in {\mathbb C}$ for $j\in J$.
Thus $b(u\otimes v)=u\otimes (\sum_{j\in J}\alpha_{j}v_{j})\in u\otimes V_{i}$.
Therefore
\begin{eqnarray}
b(u\otimes v)\in u\otimes V_{i}
\end{eqnarray}
for all $u\in U_{i}, \;v\in V_{i}$. It follows that 
$b\in \coprod_{i\in I}{\mathbb C}\otimes gl_{J}(V_{i})$. This proves (\ref{eC(A)1}).

Since $V_{i}$ $(i\in I)$ are irreducible graded $C(A)$-modules
(by Proposition \ref{p2}),
using (5.5) we have
\begin{eqnarray}
C(C(A))=\coprod_{i\in I}gl_{J}(U_{i})\otimes {\mathbb C}.
\end{eqnarray}
Then it follows from Theorem \ref{tAW} that $C(C(A))=\rho(\tilde{A})$.
\end{proof}

By Corollary \ref{cdirectsum}, $\oplus_{i\in I}({\mathbb C}\otimes
gl_{J}(V_{i}))$ is a rational graded associative algebra if $I$ is finite.
In view of Proposition \ref{pdual1}, $C(A)$ with a certain topology
is a rational graded associative algebra. Next we shall
prove that $\oplus_{i\in I}({\mathbb C}\otimes
gl_{J}(V_{i}))$ is a topological subspace of $gl_{J}(W)$,
so that $C(A)$ with the induced topology from $gl_{J}(W)$	
is rational.

Let $U'=\coprod_{n\in {\mathbb Z}}U'(n)$ be a graded
subspace of a ${\mathbb Z}$-graded vector space $U$. In Section 2, 
we have defined a Hausdorff topology on $U$.
Then we have two topologies on $U'$,
which are the induced topology from $U$ and the topology 
defined by using its own homogeneous subspaces.
These two topologies on $U'$ are in fact the same
because for $u\in U, \;k\in {\mathbb Z}$,
$$U'\cap \left(u+\sum_{n\ge k}U(n)\right)
=\cup _{u'}\left(u'+\sum_{n\ge k}U'(n)\right),$$
where $u'$ runs through the set on the left hand side.

\begin{lemma}\label{ltopolo}
Suppose that $W=W_{1}\oplus \cdots \oplus W_{r}$, where 
$W_{1},\dots, W_{r}$ are graded subspaces of $W$.
With this decomposition we identify each $gl_{J}(W_{i})$ as
a subspace of $gl_{J}(W)$.
Then the  topological space $gl_{J}(W_{i})$ is a topological subspace of
$gl_{J}(W)$.
\end{lemma}

\begin{proof} View $\End' W_{i}$ as a graded subspace of $\End' W$.
Then from the discussion right before this lemma 
$\End'W_{i}$ is a topological subspace 
of $\End' W$. Since $gl_{J}(W_{i})$ and $gl_{J}(W)$ are
the completions of $\End'W_{i}$
and $\End' W$, respectively, $gl_{J}(W_{i})$ is a 
topological subspace of $gl_{J}(W)$.
\end{proof}

Now we have:
 
\begin{corollary}\label{ccarational}
If $I$ is finite, the commutant algebra $C(A)$ with 
the induced topology is rational and 
$\Hom_{A}(U_{i},W)$ for $i\in I$
exhaust all inequivalent irreducible ${\mathbb N}$-graded
continuous $C(A)$-modules.
\end{corollary}

\begin{proof} In view of Corollary \ref{cdirectsum}, we need to prove that 
$\coprod_{i\in I}({\mathbb C}\otimes gl_{J}(V_{i}))$ 
is a topological subspace of $gl_{J}(W)$.
Since ${\mathbb C}\otimes gl_{J}(V_{i})$ is the completion of
${\mathbb C}\otimes \End' V_{i}$ and $\End'(U_{i}\otimes V_{i})$
is a topological subspace of $gl_{J}(W)$ (by Lemma \ref{ltopolo}),
it suffices to prove that
${\mathbb C}\otimes \End' V_{i}$ is a topological subspace of 
$\End'(U_{i}\otimes V_{i})$. This is true because
for $m\in {\mathbb Z}$, 
$$({\mathbb C}\otimes \End' V_{i})_{m}
=\oplus_{n}\left({\mathbb C}\otimes \Hom(V_{i}(n),V_{i}(m+n))\right)$$
with respect to this very grading is a graded subspace of 
$$(\End'W)_{m}=\oplus_{n}\Hom(W(n),W(m+n))$$
Then the proof is complete.
\end{proof}

Let $I=ann_{\tilde{A}} \rho$. Then $A$ is a subalgebra of $\tilde{A}/I$.
We cannot prove that $A$ is rational with the induced topology if 
$I$ is finite,
however, the extension $\tilde{A}/I$ of $A$ is rational.


\section{An application to vertex operator algebras}
Throughout this section, $V$ will be a vertex operator algebra.
We shall use standard definitions and notations as defined in [FLM] 
and [FHL].

In [FZ], Frenkel and Zhu associated a ${\mathbb Z}$-graded topological algebra
$U(V)$ to $V$ such that the category of $V$-modules and
the category of continuous 
$U(V)$-modules of a certain type are equivalent.
Now we recall the definition of $U(V)$.
Let $A$ be the free associative algebra with unit
on the vector space $\hat{V}=V\otimes {\mathbb C}[t,t^{-1}]$.
By assigning a degree $\wt v-n-1$ to each element $v\otimes t^{n}$ for
homogeneous $v\in V, \;n\in {\mathbb Z}$, $A$ becomes a ${\mathbb Z}$-graded
algebra. 
Since the Jacobi identity ([B], [FLM], [FHL]) involves infinite sums,
it is necessary to complete $A$ before
to enforce the Jacobi identity relation. 
Let $\tilde{A}$ be the completion defined in Section 2.
Let $Jac$ be the two-sided
ideal of $\tilde{A}$ generated by the following relations:
\begin{eqnarray}
& &{\bf 1}(n)=\delta_{n,-1},\\
& &(L(-1)u)(n)=-nu(n-1),\\
& &\sum_{i\ge 0}\left(\begin{array}{c}l\\i\end{array}\right)
(-1)^{i}u(m+l-i)v(n+i)
-\sum_{i\ge 0}\left(\begin{array}{c}l\\i\end{array}\right)
(-1)^{l+i}v(n+l-i)u(m+i)\\
& &=\sum_{i\ge 0}\left(\begin{array}{c}m\\i\end{array}\right)(u_{l+i}v)(m+n-i)
\nonumber
\end{eqnarray}
for $u,v\in V,\;l,m,n\in {\mathbb Z}$, where $v(n)=v\otimes t^{n}$.
Set $U(V)=\tilde{A}/Jac$. This is the Frenkel and Zhu's 
universal enveloping algebra associated to $V$.
Then
any $V$-module is a natural continuous $U(V)$-module where
$u(m)$ is represented by $u_{m}$ for $u\in V, \;m\in {\mathbb Z}$.


As a corollary of Theorem \ref{tAW} we have:

\begin{theorem}\label{t3}
Let $V$ be a vertex operator algebra and let $W$ be an 
irreducible $V$-module.
Let $\rho$ be the algebra homomorphism from $U(V)$ to $gl_{J}(W)$.
Then $\rho(U(V))=gl_{J}(W)$.
\end{theorem}

\begin{proof} Since $W$ is an irreducible $V$-module, $W$ is an irreducible
$A$-module. By Theorem \ref{tAW}, $\rho(\tilde{A})=gl_{J}(W)$.
Since $W$ is a $V$-module, $Jac\cdot W=0$. Thus
$$\rho(U(V))=\rho(\tilde{A})=gl_{J}(W).$$
This completes the proof. 
\end{proof}

Theorem \ref{t3} and Chinese Remainder Theorem 
immediately imply:

\begin{corollary}\label{crational} 
Let $W=W_{1}\oplus \cdots \oplus W_{r}$, where $W_{i}$ are inequivalent
irreducible $V$-modules. Then
\begin{eqnarray}
U(V)/I=gl_{J}(W_{1})\oplus \cdots \oplus gl_{J}(W_{r}),
\end{eqnarray}
where $I=\cap_{i=1}^{r}ann_{U(V)}W_{i}$. 
\end{corollary}

Recall from [Z] (cf. [DLM2]) that $V$ is said to be {\em rational}
if any continuous ${\mathbb N}$-graded $U(V)$-module is 
completely reducible. We end this paper with the following:

\begin{conjecture}\label{conj} 
Let $V$ be a rational vertex operator algebra. Then
\begin{eqnarray}
U(V)=gl_{J}(W_{1})\oplus \cdots 
\oplus gl_{J}(W_{r}),
\end{eqnarray}
where $W_{1},\dots,W_{r}$
form a complete set of representatives of equivalent classes of
irreducible $V$-modules.
\end{conjecture}

\end{document}